\documentclass{article}
\usepackage{amsfonts}
\usepackage{amsmath}

\input{tcilatex}
\begin{document}

\begin{center}
\textbf{A New Approach for Higher Order Difference Equations and }

\textbf{Eigenvalue problems via Physical Potentials}

\bigskip

Erdal BAS$^{a}$ and Ramazan OZARSLAN$^{a\ast }$

\bigskip

$^{a}$\textit{Firat University, Science Faculty, Department of Mathematics,
23119 Elazig/Turkey}

\bigskip

\textit{e-mail: erdalmat@yahoo.com, }$^{\ast }$\textit{%
ozarslanramazan@gmail.com}

\bigskip

\textbf{Abstract}
\end{center}

{\footnotesize In this study, we give the variation of parameters method
from a different viewpoint for the }${\footnotesize Nth}${\footnotesize \
order inhomogeneous linear ordinary difference equations with constant
coefficient by means of delta exponential function }$e_{p}\left( t,s\right) $%
{\footnotesize . Advantage of this new approachment is to enable us to
investigate the solution of difference equations in the closed form. Also,
the method is supported with three difference eigenvalue problems; the
second-order Sturm-Liouville problem, which is called also one dimensional
Schr\"{o}dinger equation, having Coulomb potential, hydrogen atom equation,
and the fourth-order relaxation difference equations}$.$ {\footnotesize We
find sum representation of solution for the second order discrete
Sturm-Liouville problem having Coulomb potential, hydrogen atom equation,
and analytical solution of the fourth order discrete relaxation problem by
the variation of parameters method via delta exponential and delta
trigonometric functions .}

\bigskip


\bigskip


\bigskip

\begin{center}
\textbf{1. Introduction}

\bigskip
\end{center}

Variation of parameters method, that is a general method in the solution
methods of inhomogeneous linear ordinary differential equations, firstly was
introduced by Euler and Lagrange while they were studying the celestial
bodies and orbital elements \cite{wiki,wiki2}.

Later, this method was adapted to solve inhomogeneous linear ordinary
difference equations for the following equation types in \cite{peterson,
bender}%
\begin{equation}
x_{n+N}+p_{N-1}x_{n+N-1}+\cdots +p_{0}x_{n}=q\left( n\right) .  \tag{1}
\label{i1}
\end{equation}%
Recently, delta exponential function has been described in \cite{peter}
similarly to exponential function in the continuous case. $e^{pt},$ $p\in 
\mathbb{R}
$, is the solution of the following problem 
\begin{eqnarray*}
y^{\prime }\left( t\right) &=&py\left( t\right) , \\
y\left( 0\right) &=&1,
\end{eqnarray*}%
and delta exponential function $e_{p}\left( n,s\right) $ is the solution of
the following problem%
\begin{eqnarray*}
\Delta x\left( n\right) &=&p\left( n\right) x\left( n\right) , \\
x\left( s\right) &=&1.
\end{eqnarray*}%
Also, delta exponential function $e_{p}\left( n,s\right) $ is used to find
the homogeneous solution of linear difference equations with constant
coefficient as with differential equations and besides, variation of
constants formula is given for the first order linear difference equations
in \cite{peter}.

In this study, we generalize the method for $Nth$\ order inhomogeneous
linear ordinary difference equations with constant coefficient by means of
delta exponential function $e_{p}\left( t,s\right) $ and consider the
following equation in a closed form differently from the equation $\left( %
\ref{i1}\right) $%
\begin{equation}
\Delta ^{N}x\left( n\right) +r_{N-1}\Delta ^{N-1}x\left( n\right) +\cdots
+r_{0}x_{n}=q\left( n\right) .  \tag{2}  \label{i2}
\end{equation}

We present a different approach to the variation of parameters method by
means of delta exponential function for higher order inhomogeneous linear
ordinary difference equations with constant coefficient. This approachment
enables us to investigate the solution of difference equations in the closed
form as given in $\left( \ref{i2}\right) $, otherwise equation $\left( \ref%
{i2}\right) $ has to be expanded by using binomial expansion of the
difference operator%
\begin{equation*}
\Delta ^{N}x\left( t\right) =\sum_{k=0}^{N}\left( -1\right) ^{k}\left( 
\begin{array}{c}
N \\ 
k%
\end{array}%
\right) x\left( t+N-k\right) .
\end{equation*}

In the last section, we give three examples for explaining the method, also
specifically we find the representation of solution of discrete
Sturm-Liouville problem,{\footnotesize \ }which is called also one
dimensional Schr\"{o}dinger equation, having Coulomb potential, discrete
hydrogen atom equation, and the fourth-order discrete relaxation equation.

Sturm-Liouville equations play an important role in mathematical physics.
Lately, Sturm-Liouville differential and difference equations have been
considered similarly to the continuous counterpart $\left[ 1-10,14\right] $.

Hydrogen atom equation is studied by \cite{ulusoy, bash, everitt, adkins}.
Hydrogen atom equation is used in quantum mechanics for determining energy
levels of hydrogen atom \cite{ulusoy}. Hydrogen atom equation is defined as
follows%
\begin{equation*}
y^{\prime \prime }+\left( \lambda -\dfrac{l\left( l+1\right) }{x^{2}}+\dfrac{%
2}{x}-q\left( x\right) \right) y=0.
\end{equation*}%
Then, let's introduce discrete hydrogen atom equation%
\begin{equation*}
\Delta ^{2}x\left( n-1\right) +\left( \lambda -q\left( n\right) +\frac{2}{n}-%
\frac{l\left( l+1\right) }{n^{2}}\right) x\left( n\right) =0,\text{ }%
n=1,...,b,
\end{equation*}%
where $l$ is a positive integer or zero, $v\left( n\right) \in l^{2}\left[
0,b\right] ,$ $q\left( n\right) ,$ $b,$ $\lambda $ and $n$ is as defined
above, $-q\left( n\right) +\dfrac{2}{n}-\dfrac{l\left( l+1\right) }{n^{2}}$
are called potential function.

Now, let's define briefly Sturm-Liouville operator having Coulomb potential.
Motion of electrons moving under the Coulomb potential has importance in
quantum theory. This problem is used for finding energy levels for hydrogen
atom and single valence electron atoms. Time-dependent Schr\"{o}dinger
equation is as follows 
\begin{equation*}
i\tilde{h}\frac{\partial \omega }{\partial t}=-\frac{\tilde{h}^{2}}{2m}\frac{%
\partial ^{2}\omega }{\partial x^{2}}+U\left( x,y,z\right) \omega ,\text{%
\qquad }\int \limits_{R^{3}}\left \vert \omega \right \vert ^{2}dxdydz=1,
\end{equation*}%
from here, in consequence of some transformations, we obtain Sturm-Liouville
equation having Coulomb potential%
\begin{equation*}
-y^{\prime \prime }+\left[ \frac{A}{x}+q\left( x\right) \right] y=\lambda y,
\end{equation*}%
where $\lambda $ is a parameter which corresponds to the energy \cite{bilo}.
The following problem%
\begin{eqnarray*}
-\Delta ^{2}x\left( n-1\right) +\left( \frac{A}{n}+q\left( n\right) \right)
x\left( n\right)  &=&\lambda u\left( n\right) ,\text{ }n=a,...,b, \\
x\left( a-1\right) +hx\left( a\right)  &=&0,
\end{eqnarray*}%
is called discrete Sturm-Liouville problem having Coulomb potential.

\newpage 

\begin{center}
\textbf{2.} \textbf{Preliminaries}
\end{center}

\noindent \textbf{Definition 1.} \cite{peter} Let's define regressive
functions,%
\begin{equation*}
\Re =\left \{ p:%
\mathbb{N}
_{a}\rightarrow 
\mathbb{R}
\text{ such that }1+p\left( n\right) \neq 0\text{ for }n\in 
\mathbb{N}
_{a}\right \} .
\end{equation*}

\bigskip

\noindent \textbf{Theorem 2. }\cite{peter} Let's define delta exponential
function. Suppose that $p\in \Re $ and $s\in 
\mathbb{N}
_{a}$, then%
\begin{equation*}
e_{p}\left( n,s\right) =\left \{ 
\begin{array}{cc}
\prod \limits_{\tau =s}^{n-1}\left[ 1+p\left( \tau \right) \right] , & n\in 
\mathbb{N}
_{s} \\ 
\prod \limits_{\tau =n}^{s-1}\left[ 1+p\left( \tau \right) \right] ^{-1}, & 
n\in 
\mathbb{N}
_{a}^{s-1},%
\end{array}%
\right.
\end{equation*}%
where $\prod \limits_{\tau =a}^{b}.=1,$ if $a>b.$

\bigskip

\noindent \textbf{Theorem 3. }\cite{peter} Suppose that $p,$ $q\in \Re $ and 
$n,$ $s\in 
\mathbb{N}
_{a}.$ Then

\noindent $\left( \text{\textbf{i}}\right) $ $e_{0}\left( n,s\right) =1;$

\noindent $\left( \text{\textbf{ii}}\right) $ $\Delta e_{p}\left( n,s\right)
=p\left( n\right) e_{p}\left( n,s\right) ;$

\noindent $\left( \text{\textbf{iii}}\right) $ $e_{p}\left( n,s\right)
e_{q}\left( n,s\right) =e_{p\oplus q}\left( n,s\right) ,$

\noindent where $p\oplus q=p+q+pq.$

\bigskip

\noindent \textbf{Definition 4. }\cite{peter} Let's define delta sine and
cosine functions as follows,%
\begin{equation*}
\cos _{p}\left( n,a\right) =\frac{e_{ip}\left( n,a\right) +e_{-ip}\left(
n,a\right) }{2},\text{ }\sin _{p}\left( n,a\right) =\frac{e_{ip}\left(
n,a\right) -e_{-ip}\left( n,a\right) }{2i},
\end{equation*}%
where $n\in 
\mathbb{N}
_{a},$ $\pm ip\in \Re .$

\bigskip

\noindent \textbf{Definition 5. }\cite{peter} Let's define delta integral.
Suppose $f:%
\mathbb{N}
_{a}\rightarrow 
\mathbb{R}
$ and $c\leq d$, $c,$ $d\in 
\mathbb{N}
_{a},$ then%
\begin{equation*}
\int \limits_{c}^{d+1}f\left( n\right) \Delta n=\sum_{n=c}^{d}f\left(
n\right) ,
\end{equation*}%
where $\sum \limits_{n=c}^{d}.=0,$ if $c>d.$

\bigskip

\noindent \textbf{Definition 6. }\cite{peter} If $F\left( n\right) $ is
delta integral of $f\left( n\right) $, then%
\begin{equation*}
\int \limits_{a}^{b}f\left( n\right) \Delta n=F\left( b\right) -F\left(
a\right) .
\end{equation*}

Let's consider the second order linear homogeneous ordinary difference
equation with constant coefficients as follows,%
\begin{equation}
\Delta ^{2}y\left( n\right) +p\Delta y\left( n\right) +qy\left( n\right) =0,%
\text{ }n\in 
\mathbb{N}
_{a},  \tag{3}  \label{1}
\end{equation}%
where $p,$ $q\in 
\mathbb{R}
$ hold $p\neq 1+q.$

\bigskip

\noindent \textbf{Theorem 7. }\cite{peter} Characteristic equation of $%
\left( \ref{1}\right) ,$ by the help of delta exponential function, is given
by%
\begin{equation*}
m^{2}+pm+q=0,
\end{equation*}%
let $m_{1},$ $m_{2}$ are distinct characteristic roots of the characterstic
equation , so%
\begin{equation}
y\left( n\right) =c_{1}e_{m_{1}}\left( n,a\right) +c_{2}e_{m_{2}}\left(
n,a\right) ,  \tag{4}  \label{2}
\end{equation}%
where $c_{1},$ $c_{2}$ are constants, is a general solution of $\left( \ref%
{1}\right) .$

\bigskip

\noindent \textbf{Theorem 8. }\cite{peter} Let the characteristic roots are
complex pair, $m_{1,2}=\alpha \pm i\beta ,$ $\alpha \neq -1,$ $\beta >0,$ so 
\begin{equation}
y\left( n\right) =c_{1}e_{\alpha }\left( n,a\right) \cos _{\gamma }\left(
n,a\right) +c_{2}e_{\alpha }\left( n,a\right) \sin _{\gamma }\left(
n,a\right) ,  \tag{5}  \label{3}
\end{equation}%
where $\gamma =\dfrac{\beta }{1+\alpha },$ is a general solution of $\left( %
\ref{1}\right) .$

\bigskip

\noindent \textbf{Theorem 9. }\cite{peter} Let the characteristic roots are
double roots, $m_{1}=m_{2}=r,$ so 
\begin{equation}
y\left( n\right) =c_{1}e_{r}\left( n,a\right) +c_{2}\left( n-a\right)
e_{r}\left( n,a\right) ,  \tag{6}  \label{4}
\end{equation}%
is a general solution of $\left( \ref{1}\right) .$

\begin{center}
\textbf{3.} \textbf{Main Results }
\end{center}

\noindent \textbf{3.1. Analysis of the Method}

\bigskip

Let's reconsider the equation $\left( \ref{i2}\right) ;$%
\begin{equation*}
\Delta ^{N}x\left( n\right) +r_{N-1}\Delta ^{N-1}x\left( n\right) +\cdots
+r_{0}x_{n}=q\left( n\right) .
\end{equation*}%
If we change the variable $x\left( n\right) =e_{m}\left( n,0\right) $ and
consider the homogeneous part of $\left( \ref{i2}\right) ,$ then we have the
following characteristic equation%
\begin{equation*}
m^{N}+r_{N-1}m^{N-1}+\cdots +r_{0}=0,
\end{equation*}%
and let its roots are $m_{1},m_{2},...,m_{N}.$ Hence, we have the
homogeneous solution as follows,%
\begin{equation*}
x\left( n\right) =c_{1}x_{1}\left( n\right) +c_{2}x_{2}\left( n\right)
+\cdots +c_{N}x_{N}\left( n\right) ,
\end{equation*}%
where $x_{1}\left( n\right) =e_{m_{1}}\left( n,0\right) ,x_{2}\left(
n\right) =e_{m_{2}}\left( n,0\right) ,...,x_{N}\left( n\right)
=e_{m_{N}}\left( n,0\right) $ is a linearly independent set of solutions.

From here, let's take a set of new solution functions for the variation of
parameters method,\linebreak $v_{1}\left( n\right) ,v_{2}\left( n\right)
,...,v_{N}\left( n\right) $ and so, let's assume that following equation 
\begin{equation*}
X\left( n\right) =v_{1}\left( n\right) x_{1}\left( n\right) +v_{2}\left(
n\right) x_{2}\left( n\right) +\cdots +v_{N}\left( n\right) x_{N}\left(
n\right)
\end{equation*}%
is a solution of nonhomogeneous part of $\left( \ref{i2}\right) .$

For finding the parameters $v_{1}\left( n\right) ,v_{2}\left( n\right)
,...,v_{N}\left( n\right) ,$ firstly let's take 
\begin{eqnarray*}
\Delta X\left( n\right) &=&\left[ \Delta v_{1}\left( n\right) x_{1}\left(
n+1\right) +\Delta v_{2}\left( n\right) x_{2}\left( n+1\right) +\cdots
+\Delta v_{N}\left( n\right) x_{N}\left( n+1\right) \right] + \\
&&\left[ v_{1}\left( n\right) \Delta x_{1}\left( n\right) +v_{2}\left(
n\right) \Delta x_{2}\left( n\right) +\cdots +v_{N}\left( n\right) \Delta
x_{N}\left( n\right) \right] ,
\end{eqnarray*}%
let's assume that the first bracketed part at the right hand side of the
equation above equals to zero,%
\begin{equation*}
\Delta v_{1}\left( n\right) x_{1}\left( n+1\right) +\Delta v_{2}\left(
n\right) x_{2}\left( n+1\right) +\cdots +\Delta v_{N}\left( n\right)
x_{N}\left( n+1\right) =0,
\end{equation*}%
and so we have,%
\begin{equation*}
\Delta X\left( n\right) =v_{1}\left( n\right) \Delta x_{1}\left( n\right)
+v_{2}\left( n\right) \Delta x_{2}\left( n\right) +\cdots +v_{N}\left(
n\right) \Delta x_{N}\left( n\right) ,
\end{equation*}%
then, let's take the difference of the equality above%
\begin{eqnarray*}
\Delta ^{2}X\left( n\right) &=&\left[ \Delta v_{1}\left( n\right) \Delta
x_{1}\left( n+1\right) +\Delta v_{2}\left( n\right) \Delta x_{2}\left(
n+1\right) +\cdots +\Delta v_{N}\left( n\right) \Delta x_{N}\left(
n+1\right) \right] + \\
&&\left[ v_{1}\left( n\right) \Delta ^{2}x_{1}\left( n\right) +v_{2}\left(
n\right) \Delta ^{2}x_{2}\left( n\right) +\cdots +v_{N}\left( n\right)
\Delta ^{2}x_{N}\left( n\right) \right] ,
\end{eqnarray*}%
let's assume that%
\begin{equation*}
\Delta v_{1}\left( n\right) \Delta x_{1}\left( n+1\right) +\Delta
v_{2}\left( n\right) \Delta x_{2}\left( n+1\right) +\cdots +\Delta
v_{N}\left( n\right) \Delta x_{N}\left( n+1\right) =0,
\end{equation*}%
so we have%
\begin{equation*}
\Delta ^{2}X\left( n\right) =v_{1}\left( n\right) \Delta ^{2}x_{1}\left(
n\right) +v_{2}\left( n\right) \Delta ^{2}x_{2}\left( n\right) +\cdots
+v_{N}\left( n\right) \Delta ^{2}x_{N}\left( n\right) ,
\end{equation*}%
proceeding in this fashion, we have%
\begin{equation*}
\Delta ^{N}X\left( n\right) =v_{1}\left( n\right) \Delta ^{N}x_{1}\left(
n\right) +v_{2}\left( n\right) \Delta ^{N}x_{2}\left( n\right) +\cdots
+v_{N}\left( n\right) \Delta ^{N}x_{N}\left( n\right) .
\end{equation*}%
If we formulate our assumptions above, we have%
\begin{equation*}
\Delta ^{k}X\left( n\right) =v_{1}\left( n\right) \Delta ^{k}x_{1}\left(
n\right) +v_{2}\left( n\right) \Delta ^{k}x_{2}\left( n\right) +\cdots
+v_{N}\left( n\right) \Delta ^{k}x_{N}\left( n\right) ,\text{ }k=1,2,...,N-1,
\end{equation*}%
\begin{equation*}
\Delta ^{k}v_{1}\left( n\right) \Delta x_{1}\left( n+1\right) +\Delta
^{k}v_{2}\left( n\right) \Delta x_{2}\left( n+1\right) +\cdots +\Delta
^{k}v_{N}\left( n\right) \Delta x_{N}\left( n+1\right) =0,\text{ }%
k=0,1,...,N-2,
\end{equation*}%
if obtained equalities above is written in $\left( \ref{i2}\right) ,$ then
we have the following equation system 
\begin{eqnarray*}
\Delta v_{1}\left( n\right) x_{1}\left( n+1\right) +\Delta v_{2}\left(
n\right) x_{2}\left( n+1\right) +\cdots +\Delta v_{N}\left( n\right)
x_{N}\left( n+1\right) &=&0, \\
\Delta v_{1}\left( n\right) \Delta x_{1}\left( n+1\right) +\Delta
v_{2}\left( n\right) \Delta x_{2}\left( n+1\right) +\cdots +\Delta
v_{N}\left( n\right) \Delta x_{N}\left( n+1\right) &=&0, \\
&&\vdots \\
\Delta v_{1}\left( n\right) \Delta ^{N-2}x_{1}\left( n+1\right) +\Delta
v_{2}\left( n\right) \Delta ^{N-2}x_{2}\left( n+1\right) +\cdots +\Delta
v_{N}\left( n\right) \Delta ^{N-2}x_{N}\left( n+1\right) &=&0, \\
\Delta v_{1}\left( n\right) \Delta ^{N-1}x_{1}\left( n+1\right) +\Delta
v_{2}\left( n\right) \Delta ^{N-1}x_{2}\left( n+1\right) +\cdots +\Delta
v_{N}\left( n\right) \Delta ^{N-1}x_{N}\left( n+1\right) &=&q\left( n\right)
.
\end{eqnarray*}%
If we solve this system by Cramer rule, we have the Casoratian,%
\begin{equation*}
\det \left \vert 
\begin{array}{cccc}
x_{1}\left( n+1\right) & x_{2}\left( n+1\right) & \cdots & x_{N}\left(
n+1\right) \\ 
\Delta x_{1}\left( n+1\right) & \Delta x_{2}\left( n+1\right) & \cdots & 
\Delta x_{N}\left( n+1\right) \\ 
\vdots & \vdots & \ddots & \vdots \\ 
\Delta ^{N-1}x_{1}\left( n+1\right) & \Delta ^{N-1}x_{2}\left( n+1\right) & 
\cdots & \Delta ^{N-1}x_{N}\left( n+1\right)%
\end{array}%
\right \vert =W\left( x_{1},x_{2},\cdots ,x_{N}\right) \left( n+1\right)
\end{equation*}%
and assuming the Casoratian is different from zero. Let $W_{i}$ correspond
the determinant of the $i$th column of the Casoratian with the column $%
\left( 0,0,0,...,0,1\right) $ and so, solution of the system as follows%
\begin{equation*}
\Delta v_{1}\left( n\right) =\frac{q\left( n\right) W_{1}\left( n\right) }{%
W\left( n+1\right) },\Delta v_{2}\left( n\right) =\frac{q\left( n\right)
W_{2}\left( n\right) }{W\left( n+1\right) },\cdots ,\Delta v_{N}\left(
n\right) =\frac{q\left( n\right) W_{N}\left( n\right) }{W\left( n+1\right) },
\end{equation*}%
from here we have the parameters as follows%
\begin{equation*}
v_{1}\left( n\right) =\int \frac{q\left( n\right) W_{1}\left( n\right) }{%
W\left( n+1\right) }\Delta n,v_{2}\left( n\right) =\int \frac{q\left(
n\right) W_{2}\left( n\right) }{W\left( n+1\right) }\Delta n,\cdots
,v_{N}\left( n\right) =\int \frac{q\left( n\right) W_{N}\left( n\right) }{%
W\left( n+1\right) }\Delta n.
\end{equation*}%
Finally, the particular solution is as follows%
\begin{equation*}
X\left( n\right) =x_{1}\left( n\right) \int \frac{q\left( n\right)
W_{1}\left( n\right) }{W\left( n+1\right) }\Delta n+x_{2}\left( n\right)
\int \frac{q\left( n\right) W_{2}\left( n\right) }{W\left( n+1\right) }%
\Delta n+\cdots +x_{N}\left( n\right) \int \frac{q\left( n\right)
W_{N}\left( n\right) }{W\left( n+1\right) }\Delta n.
\end{equation*}

\newpage 

\noindent \textbf{3.2. Numerical Results and Discussions of Some Discrete
Eigenvalue Problems Having Physical Potential}

In this section, a new version of the variation of parameters method is
applied by using delta exponential function. First of all, we consider the
second-order Sturm-Liouville problem, which is called also one dimensional
Schr\"{o}dinger equation, having Coulomb potential, hydrogen atom equation,
and the fourth-order relaxation difference equations$.$ We find sum
representation of solution for the second order discrete Sturm-Liouville
problem having Coulomb potential, hydrogen atom equation, and analytical
solution of the fourth order discrete relaxation problem by means of
variation of parameters method by using delta exponential function.\textbf{\ 
}

\bigskip 

\noindent \textbf{3.2.1. Discrete Hydrogen Atom Equation and Discrete
Sturm-Liouville Equation Having Coulomb Potential}

First of all, let's consider the following second order hydrogen atom
equation,%
\begin{equation}
-\Delta ^{2}x\left( n-1\right) +\left( \frac{A}{n}+q\left( n\right) \right)
x\left( n\right) =\lambda x\left( n\right) ,  \tag{7}  \label{a2}
\end{equation}%
with the boundary conditions, 
\begin{equation}
x\left( 0\right) =x\left( b\right) =0,  \tag{8}  \label{a3}
\end{equation}%
has a unique solution $x\left( n\right) $ as follows,$n\in \left[ 0,b\right]
,$ $n$ is a finite integer, $x\left( n\right) \in l^{2}\left[ 0.b\right] $%
\begin{equation*}
x\left( n,\lambda \right) =\sin n\theta +\frac{1}{\sin \theta }\int
\limits_{1}^{n+1}\left( \frac{A}{s}+q\left( s\right) \right) x\left(
s\right) \sin \left( n-s\right) \theta \Delta s.
\end{equation*}%
Secondly, let's consider the following second order Sturm-Liouville
difference equation having Coulomb potential,%
\begin{equation}
-\Delta ^{2}x\left( n-1\right) +\left( -q\left( n\right) +\frac{2}{n}-\frac{%
l\left( l+1\right) }{n^{2}}\right) x\left( n\right) =\lambda x\left(
n\right) ,  \tag{9}  \label{a4}
\end{equation}%
with the boundary conditions $\left( \ref{a3}\right) ,$ has a unique
solution $x\left( n\right) $ as follows, $n\in \left[ 0,b\right] ,$ $n$ is a
finite integer, $x\left( n\right) \in l^{2}\left[ 0.b\right] $%
\begin{equation}
x\left( n,\lambda \right) =\sin n\theta +\frac{1}{\sin \theta }\int
\limits_{1}^{n+1}\left( -q\left( s\right) +\frac{2}{s}-\frac{l\left(
l+1\right) }{s^{2}}\right) x\left( s\right) \sin \left( n-s\right) \theta
\Delta s.  \tag{10}  \label{7}
\end{equation}%
\begin{equation*}
\underset{Fig1:\text{ Eigenfunctions for the problem }\left( 7\right)
-\left( 8\right) }{\FRAME{itbpF}{3.0441in}{1.9294in}{0in}{}{}{Figure}{%
\special{language "Scientific Word";type "GRAPHIC";maintain-aspect-ratio
TRUE;display "USEDEF";valid_file "T";width 3.0441in;height 1.9294in;depth
0in;original-width 3in;original-height 1.8913in;cropleft "0";croptop
"1";cropright "1";cropbottom "0";tempfilename
'P32BRJ0B.wmf';tempfile-properties "XPR";}}\qquad }\underset{Fig2:\text{
Eigenfunctions for the problem }\left( 9\right) -\left( 8\right) }{\FRAME{%
itbpF}{3.2007in}{1.9225in}{0in}{}{}{Figure}{\special{language "Scientific
Word";type "GRAPHIC";maintain-aspect-ratio TRUE;display "USEDEF";valid_file
"T";width 3.2007in;height 1.9225in;depth 0in;original-width
3in;original-height 1.7919in;cropleft "0";croptop "1";cropright
"1";cropbottom "0";tempfilename 'P329W508.wmf';tempfile-properties "XPR";}}}
\end{equation*}%
\begin{equation*}
\underset{Fig3:\text{Comparison of the eigenfunctions for the problems }%
\left( 7\right) -\left( 8\right) \text{ and }\left( 9\right) -\left(
8\right) }{\FRAME{itbpF}{3.0441in}{1.9294in}{0in}{}{}{Figure}{\special%
{language "Scientific Word";type "GRAPHIC";maintain-aspect-ratio
TRUE;display "USEDEF";valid_file "T";width 3.0441in;height 1.9294in;depth
0in;original-width 3in;original-height 1.8913in;cropleft "0";croptop
"1";cropright "1";cropbottom "0";tempfilename
'P329W50A.wmf';tempfile-properties "XPR";}}}
\end{equation*}

$\underset{Table1:\text{ Eigenfunctions correspond to the eigenvalue }%
\lambda =1,\text{ }}{%
\begin{tabular}{|c|c|c|}
$n$ & $x\left( n\right) $ for Hydrogen & $x\left( n\right) $ for Coulomb \\ 
$1$ & $0.866025$ & $0.866025$ \\ 
$2$ & $2.59808$ & $5.19615$ \\ 
$3$ & $4.86821$ & $10.6024$ \\ 
$4$ & $6.70353$ & $11.5276$ \\ 
$5$ & $6.86296$ & $5.24802$ \\ 
$6$ & $4.60124$ & $-4.77228$ \\ 
$13$ & $6.18061$ & $-7.29527$ \\ 
$14$ & $4.53803$ & $-11.2966$ \\ 
$15$ & $-0.105595$ & $-5.75247$ \\ 
$16$ & $-4.67793$ & $4.67243$ \\ 
$21$ & $5.26773$ & $-9.75089$ \\ 
$23$ & $-3.99289$ & $8.78504$ \\ 
$24$ & $-6.01193$ & $10.8038$ \\ 
$25$ & $-3.49672$ & $3.43627$%
\end{tabular}%
}\underset{Table2:\text{ Eigenfunctions correspond to the eigenvalue }%
\lambda =2-\sqrt{2},\text{ }}{\qquad 
\begin{tabular}{|c|c|c|}
$n$ & $x\left( n\right) $ for Hydrogen & $x\left( n\right) $ for Coulomb \\ 
$1$ & $0.707107$ & $0.707107$ \\ 
$2$ & $2.41421$ & $4.53553$ \\ 
$3$ & $5.62132$ & $11.182$ \\ 
$4$ & $10.6548$ & $17.7341$ \\ 
$5$ & $17.4379$ & $20.5481$ \\ 
$6$ & $25.2922$ & $17.227$ \\ 
$13$ & $-7.27102$ & $6.52242$ \\ 
$14$ & $-22.8474$ & $16.3393$ \\ 
$15$ & $-32.7782$ & $19.1177$ \\ 
$16$ & $-34.1566$ & $13.5942$ \\ 
$21$ & $32.7366$ & $-17.9177$ \\ 
$23$ & $19.0974$ & $2.89076$ \\ 
$24$ & $0.835118$ & $19.0276$ \\ 
$25$ & $-17.7111$ & $15.2384$%
\end{tabular}%
}$

\noindent $\underset{Fig4:\text{ Comparison of datas in Table1}}{\FRAME{itbpF%
}{3.346in}{2.2874in}{0in}{}{}{Figure}{\special{language "Scientific
Word";type "GRAPHIC";maintain-aspect-ratio TRUE;display "USEDEF";valid_file
"T";width 3.346in;height 2.2874in;depth 0in;original-width
3in;original-height 2.0418in;cropleft "0";croptop "1";cropright
"1";cropbottom "0";tempfilename 'P329QX04.wmf';tempfile-properties "XPR";}}}%
\underset{Fig5:\text{ Comparison of datas in Table2}}{\FRAME{itbpF}{3.3572in%
}{2.2857in}{0in}{}{}{Figure}{\special{language "Scientific Word";type
"GRAPHIC";maintain-aspect-ratio TRUE;display "USEDEF";valid_file "T";width
3.3572in;height 2.2857in;depth 0in;original-width 3in;original-height
2.034in;cropleft "0";croptop "1";cropright "1";cropbottom "0";tempfilename
'P329W506.wmf';tempfile-properties "XPR";}}}$

Suppose that $q\left( n\right) =\frac{1}{\sqrt{n}},$ $A=1,$ $l=2$ in the
figures and tables above.

\newpage 

\textbf{Proof. }Firstly, we study to find the general solution of the
equation $\left( \ref{a2}\right) $ by the variation of parameters method
without using the boundary conditions. Homogenous part of $\left( \ref{a2}%
\right) $ is as follows,%
\begin{equation*}
\Delta ^{2}x\left( n-1\right) +\lambda x\left( n\right) =0.
\end{equation*}%
By using delta exponential function, we have the characteristic equation, 
\begin{equation*}
\frac{m^{2}}{1+m}+\lambda =0,
\end{equation*}%
and characteristic roots are as follows,%
\begin{equation*}
m_{1,2}=\frac{-\lambda \pm \sqrt{\lambda \left( \lambda -4\right) }}{2},
\end{equation*}%
where $m_{1,2}\in \Re .$ So, the homogeneous solution is found from $\left( %
\ref{2}\right) $ as follows, 
\begin{equation}
x_{h}\left( n\right) =c_{1}e_{m_{1}}\left( n,0\right) +c_{2}e_{m_{2}}\left(
n,0\right) .  \tag{11}  \label{8}
\end{equation}%
Then, if we apply the variation of parameters method, lineary independent
solutions is found as $e_{m_{1}}\left( n,0\right) ,$ $e_{m_{2}}\left(
n,0\right) .$

From here, if we take the constants as parameters $v_{1}\left( n\right) $
and $v_{2}\left( n\right) ,$ then we find the variables by Cramer rule, let $%
q_{1}\left( n\right) =\frac{A}{n}+q\left( n\right) $%
\begin{equation*}
\Delta v_{1}\left( n-1\right) =\frac{%
\begin{vmatrix}
-q_{1}\left( n\right) x\left( n\right) & \Delta e_{m_{2}}\left( n,0\right)
\\ 
0 & e_{m_{2}}\left( n,0\right)%
\end{vmatrix}%
}{W\left( e_{m_{1}}\left( n,0\right) ,e_{m_{2}}\left( n,0\right) \right) },
\end{equation*}%
where $W$ is Casoratian,%
\begin{eqnarray*}
W\left( e_{m_{1}}\left( n,0\right) ,e_{m_{2}}\left( n,0\right) \right)
&=&e_{m_{1}\oplus m_{2}}\left( n,0\right) \left( m_{2}-m_{1}\right) , \\
&=&-\sqrt{\lambda \left( \lambda -4\right) }.
\end{eqnarray*}%
Hence,%
\begin{equation*}
v_{1}\left( n\right) =\frac{1}{\sqrt{\lambda \left( \lambda -4\right) }}\int
\limits_{0}^{n+1}q_{1}\left( s\right) x\left( s\right) e_{m_{2}}\left(
s,0\right) \Delta s.
\end{equation*}%
Similarly,%
\begin{equation*}
v_{2}\left( n\right) =-\frac{1}{\sqrt{\lambda \left( \lambda -4\right) }}%
\int \limits_{0}^{n+1}q_{1}\left( s\right) x\left( s\right) e_{m_{1}}\left(
s,0\right) \Delta s.
\end{equation*}%
Finally, the general solution is found by%
\begin{equation*}
x\left( n\right) =c_{1}e_{m_{1}}\left( n,0\right) +c_{2}e_{m_{2}}\left(
n,0\right) +\frac{1}{\sqrt{\lambda \left( \lambda -4\right) }}\int
\limits_{0}^{n+1}q_{1}\left( s\right) x\left( s\right) \left[
e_{m_{1}}\left( s,0\right) e_{m_{2}}\left( n,0\right) -e_{m_{2}}\left(
s,0\right) e_{m_{1}}\left( n,0\right) \right] \Delta s.
\end{equation*}%
Now, let's continue to the proof by using the boundary conditions $\left( %
\ref{a3}\right) .$ Homogeneous solution is as in the equality $\left( \ref{8}%
\right) $. For finding untrivial solution, we analyze the eigenvalue $%
\lambda $ in four cases, these are

\noindent $i)$ $\lambda =0,$

\noindent $ii)$ $\lambda =4,$

\noindent $iii)$ $\lambda >0$ and $\lambda <4,$

\noindent $iv)$ $0<\lambda <4.$

We have trivial solutions for the first two cases, we arrive at a
contradiction for the third case and finally, if $0<\lambda <4,$ then we
have untrivial solution. The characteristic roots are complex pair and
taking $\lambda =2-2\cos \theta ,$ then we have 
\begin{equation*}
m_{1,2}=\left( -1+\cos \theta \right) \pm i\sin \theta
\end{equation*}%
So, homogeneous solution is as follows by $\left( \ref{3}\right) $%
\begin{equation*}
x\left( n\right) =c_{1}e_{\alpha }\left( n,0\right) \cos _{\gamma }\left(
n,0\right) +c_{2}e_{\alpha }\left( n,0\right) \sin _{\gamma }\left(
n,0\right) ,
\end{equation*}%
where $\alpha =-1+\cos \theta ,$ $\beta =\sin \theta ,$ $\gamma =\tan \theta
.$ If we insert $\alpha ,$ $\beta ,$ $\gamma $ in the equality above and use
Theorem 2 and Definition 4, then we have the homogenous solution of $\left( %
\ref{a2}\right) $%
\begin{equation*}
x\left( n\right) =c_{1}\cos n\theta +c_{2}\sin n\theta .
\end{equation*}%
Then, it is easily found the general solution $\left( \ref{7}\right) $ by
the variation of parameters method.

Similarly, representation of solution is obtained for $\left( \ref{a4}%
\right) $

\bigskip

\noindent \textbf{3.2.2. Fourth Order Relaxation Difference Equation}

Let's consider the following the fourth order relaxation difference equation,%
\textbf{\ }%
\begin{equation}
\Delta ^{4}x\left( n\right) -\lambda x\left( n\right) =q\left( n\right) , 
\tag{12}  \label{t1}
\end{equation}%
with the initial conditions%
\begin{equation}
x\left( 0\right) =x\left( 1\right) =x\left( 2\right) =x\left( 3\right) =1. 
\tag{13}  \label{t2}
\end{equation}%
The problem $\left( \ref{t1}\right) -\left( \ref{t2}\right) $ has a unique
solution as follows, $n\in \left[ 0,b\right] ,$ $n$ is a finite integer, $%
x\left( n\right) \in l^{2}\left[ 0.b\right] $

\begin{center}
\noindent $%
\begin{array}{c}
x\left( n\right) 
\begin{array}{c}
=%
\end{array}%
\dfrac{\left( (1-is)^{n}+(1+is)^{n}\right) \left(
s^{6}+q_{0}-s^{2}(1+q_{0})-2q_{1}+q_{2}\right) }{4s^{2}\left(
-1+s^{4}\right) }+%
\end{array}%
$ 
\begin{eqnarray}
&&\frac{1}{4s^{3}\left( -1+s^{4}\right) }(1+s)^{n}\left(
s^{7}+q_{0}-s^{3}(1+q_{0})+s^{2}(q_{0}-q_{1})-3q_{1}+3q_{2}-s(q_{0}-2q_{1}+q_{2})-q_{3}\right) -
\notag \\
&&\frac{i\left( -(1-is)^{n}+(1+is)^{n}\right) \left( \left( -1+s^{2}\right)
q_{0}-\left( -3+s^{2}\right) q_{1}-3q_{2}+q_{3}\right) }{4s^{3}\left(
-1+s^{4}\right) }+  \notag \\
&&\tfrac{1}{4s^{3}\left( -1+s^{4}\right) }(1-s)^{n}\left(
s^{7}-q_{0}-s^{3}(1+q_{0})+3q_{1}+s^{2}(-q_{0}+q_{1})-3q_{2}-s(q_{0}-2q_{1}+q_{2})+q_{3}\right) -
\TCItag{14}  \label{t3} \\
&&\tfrac{1}{2}i\left( \text{-}(1\text{-}is)^{n}\text{+}(1\text{+}%
is)^{n}\right) \tsum \limits_{i=0}^{n}\left( \text{-}\tfrac{1}{4s^{3}}\left(
(1\text{-}is)^{i}+(1\text{+}is)^{i}\text{-}i\left( (1\text{-}is)^{i}\text{-}%
(1\text{+}is)^{i}\right) s\right) \left( 1\text{-}s^{2}\right) ^{1+i}\left( 1%
\text{-}s^{4}\right) ^{-1-i}q_{i}\right) \text{+}  \notag \\
&&\tfrac{1}{2}\left( (1\text{-}is)^{n}\text{+}(1\text{+}is)^{n}\right)
\tsum_{i=0}^{n}\left( \tfrac{1}{4s^{3}}\left( i\left( (1\text{-}is)^{i}\text{%
-}(1\text{+}is)^{i}\right) \text{+}\left( (1\text{-}is)^{i}\text{+}(1\text{+}%
is)^{i}\right) s\right) \left( 1\text{-}s^{2}\right) ^{1+i}\left( 1\text{-}%
s^{4}\right) ^{-1-i}q_{i}\right) \text{+}  \notag \\
&&(1+s)^{n}\sum \limits_{i=0}^{n}\left( -\frac{(1-s)^{i}\left( 1+s^{2}\right)
^{i}\left( -1+s-s^{2}+s^{3}\right) \left( 1-s^{4}\right) ^{-1-i}q_{i}}{4s^{3}%
}\right) +  \notag \\
&&(1-s)^{n}\sum \limits_{i=0}^{n}\left( -\frac{(1+s)^{i}\left( 1+s^{2}\right)
^{i}\left( 1+s+s^{2}+s^{3}\right) \left( 1-s^{4}\right) ^{-1-i}q_{i}}{4s^{3}}%
\right) .  \notag
\end{eqnarray}%
$\underset{Fig6:\text{ Eigenfunctions for the problem }\left( 11\right)
-\left( 12\right) ,\text{ }q\left( n\right) =n}{\FRAME{itbpF}{3.0441in}{%
1.9969in}{0in}{}{}{Figure}{\special{language "Scientific Word";type
"GRAPHIC";maintain-aspect-ratio TRUE;display "USEDEF";valid_file "T";width
3.0441in;height 1.9969in;depth 0in;original-width 3in;original-height
1.9579in;cropleft "0";croptop "1";cropright "1";cropbottom "0";tempfilename
'P33QOQ00.wmf';tempfile-properties "XPR";}}}$
\end{center}

\noindent $\underset{Table3:\text{ Eigenfunctions correspond to the
eigenvalue }\lambda =0.0625\qquad }{%
\begin{tabular}{|c|c|c|}
$n$ & $x\left( n\right) ,q\left( n\right) =\frac{1}{n+1}$ & $x\left(
n\right) ,q\left( n\right) =\frac{1}{\sqrt{n+1}}$ \\ 
$0$ & $1$ & $1$ \\ 
$1$ & $1$ & $1$ \\ 
$2$ & $1$ & $1$ \\ 
$3$ & $1$ & $1$ \\ 
$4$ & $16.9467$ & $2.0407$ \\ 
$5$ & $80.7556$ & $5.9615$ \\ 
$6$ & $240.321$ & $15.2739$ \\ 
$7$ & $559.519$ & $33.0829$ \\ 
$8$ & $1373.36$ & $63.0968$ \\ 
$9$ & $4308.94$ & $109.83$ \\ 
$10$ & $14838.5$ & $179.1$ \\ 
$11$ & $47386.4$ & $279.4$ \\ 
$12$ & $138351$ & $932.08$ \\ 
$13$ & $393074$ & $1377.6$%
\end{tabular}%
}\underset{Table4:\text{ Eigenfunctions correspond to the eigenvalue }%
\lambda =0.1296}{%
\begin{tabular}{|c|c|c|}
$n$ & $x\left( n\right) ,q\left( n\right) =\frac{1}{n+1}$ & $x\left(
n\right) ,q\left( n\right) =\frac{1}{\sqrt{n+1}}$ \\ 
$0$ & $1$ & $1$ \\ 
$1$ & $1$ & $1$ \\ 
$2$ & $1$ & $1$ \\ 
$3$ & $1$ & $1$ \\ 
$4$ & $5.8606$ & $2.1832$ \\ 
$5$ & $14.6518$ & $6.6554$ \\ 
$6$ & $31.0511$ & $17.3085$ \\ 
$7$ & $58.1807$ & $37.7361$ \\ 
$8$ & $100.113$ & $72.3270$ \\ 
$9$ & $162.983$ & $126.80$ \\ 
$10$ & $257.093$ & $209.55$ \\ 
$11$ & $400.414$ & $820.25$ \\ 
$12$ & $624.004$ & $1289.4$ \\ 
$13$ & $980.151$ & $2043.5$%
\end{tabular}%
}$

$\underset{}{\underset{Fig7:\text{ Comparison of datas in Table3}}{\FRAME{%
itbpF}{3.0441in}{1.9043in}{0in}{}{}{Figure}{\special{language "Scientific
Word";type "GRAPHIC";maintain-aspect-ratio TRUE;display "USEDEF";valid_file
"T";width 3.0441in;height 1.9043in;depth 0in;original-width
3in;original-height 1.8663in;cropleft "0";croptop "1";cropright
"1";cropbottom "0";tempfilename 'P2YLDM06.wmf';tempfile-properties "XPR";}}}%
\underset{Fig8:\text{ Comparison of datas in Table4}}{\FRAME{itbpF}{3.0441in%
}{1.9043in}{0in}{}{}{Figure}{\special{language "Scientific Word";type
"GRAPHIC";maintain-aspect-ratio TRUE;display "USEDEF";valid_file "T";width
3.0441in;height 1.9043in;depth 0in;original-width 3in;original-height
1.8663in;cropleft "0";croptop "1";cropright "1";cropbottom "0";tempfilename
'P2YLG808.wmf';tempfile-properties "XPR";}}}}$

\textbf{Proof. }Firstly, we study to find the general solution of the
equation $\left( \ref{t1}\right) $ by the variation of parameters method
without using the boundary conditions. Homogenous part of $\left( \ref{t1}%
\right) $ is as follows,%
\begin{equation*}
\Delta ^{4}x\left( n\right) -\lambda x\left( n\right) =0.
\end{equation*}%
By using delta exponential function, we have the characteristic equation, 
\begin{equation*}
m^{4}-\lambda =0,
\end{equation*}%
and characteristic roots are as follows,%
\begin{equation*}
m_{1,2}=\pm \sqrt[4]{\lambda },m_{3,4}=\pm i\sqrt[4]{\lambda }
\end{equation*}%
where $m_{1,2,3,4}\in \Re .$ So, the homogeneous solution is found from $%
\left( \ref{2}\right) $ as follows, 
\begin{equation*}
x_{h}\left( n\right) =c_{1}e_{m_{1}}\left( n,0\right) +c_{2}e_{m_{2}}\left(
n,0\right) +c_{3}e_{m_{3}}\left( n,0\right) +c_{4}e_{m_{4}}\left( n,0\right)
.
\end{equation*}%
Suppose that%
\begin{equation*}
\lambda =s^{4},
\end{equation*}%
and so, homogeneous solution is as follows by $\left( \ref{3}\right) $%
\begin{equation*}
x\left( n\right) =c_{1}e_{s}\left( n,0\right) +c_{2}e_{-s}\left( n,0\right)
+c_{3}e_{\alpha }\left( n,0\right) \cos _{\gamma }\left( n,0\right)
+c_{4}e_{\alpha }\left( n,0\right) \sin _{\gamma }\left( n,0\right) ,
\end{equation*}%
where $\alpha =0,$ $\beta =s,$ $\gamma =s.$ If we insert $\alpha ,$ $\beta ,$
$\gamma $ in the equality above and using Theorem 2 and Definition 4, then
we have the homogenous solution of $\left( \ref{t1}\right) $%
\begin{equation*}
x\left( n\right) =c_{1}(1-s)^{n}+c_{2}(1+s)^{n}+c_{3}\frac{%
(1+is)^{n}+(1-is)^{n}}{2}+c_{4}\frac{(1+is)^{n}-(1-is)^{n}}{2i}.
\end{equation*}%
Then, if we apply the variation of parameters method and if we take the
constants as parameters $v_{1}\left( n\right) $, $v_{2}\left( n\right) ,$ $%
v_{3}\left( n\right) $, $v_{4}\left( n\right) ,$ then we find the variables
by Cramer rule,%
\begin{eqnarray*}
v_{1} &=&\sum_{i=0}^{n}\left( \frac{-2s^{3}(1+s)^{i}\left( 1+s^{2}\right)
^{i}\left( 1+s+s^{2}+s^{3}\right) q_{i}}{8s^{6}\left( 1-s^{4}\right) ^{1+i}}%
\right) , \\
v_{2} &=&\sum_{i=0}^{n}\left( \frac{-2(1-s)^{i}s^{3}\left( 1+s^{2}\right)
^{i}\left( -1+s-s^{2}+s^{3}\right) q_{i}}{8s^{6}\left( 1-s^{4}\right) ^{1+i}}%
\right) , \\
v_{3} &=&\sum_{i=0}^{n}\left( \frac{2s^{3}\left( i\left(
(1-is)^{i}-(1+is)^{i}\right) +\left( (1-is)^{i}+(1+is)^{i}\right) s\right)
\left( 1-s^{2}\right) ^{1+i}q_{i}}{8s^{6}\left( 1-s^{4}\right) ^{1+i}}%
\right) , \\
v_{4} &=&\sum_{i=0}^{n}\left( \frac{-2s^{3}\left(
(1-is)^{i}+(1+is)^{i}-i\left( (1-is)^{i}-(1+is)^{i}\right) s\right) \left(
1-s^{2}\right) ^{1+i}q_{i}}{8s^{6}\left( 1-s^{4}\right) ^{1+i}}\right) .
\end{eqnarray*}%
Applying the initial conditions $\left( \ref{t2}\right) $, we have the
constants as follows%
\begin{eqnarray*}
c_{1} &=&\frac{1}{4s^{3}\left( -1+s^{4}\right) }\left(
s^{7}-q_{0}-s^{3}(1+q_{0})+3q_{1}+s^{2}(-q_{0}+q_{1})-3q_{2}-s(q_{0}-2q_{1}+q_{2})+q_{3}\right) ,
\\
c_{2} &=&\frac{1}{4s^{3}\left( -1+s^{4}\right) }\left(
s^{7}+q_{0}-s^{3}(1+q_{0})+s^{2}(q_{0}-q_{1})-3q_{1}+3q_{2}-s(q_{0}-2q_{1}+q_{2})-q_{3}\right) ,
\\
c_{3} &=&\frac{s^{6}+q_{0}-s^{2}(1+q_{0})-2q_{1}+q_{2}}{2s^{2}\left(
-1+s^{4}\right) }, \\
c_{4} &=&\frac{\left( -1+s^{2}\right) q_{0}-\left( -3+s^{2}\right)
q_{1}-3q_{2}+q_{3}}{2s^{3}\left( -1+s^{4}\right) }.
\end{eqnarray*}%
Hence, we obtain the general solution $\left( \ref{t3}\right) .$

\bigskip

\noindent \textbf{Conclusion}

Consequently, the variation of parameters method for higher order linear
ordinary difference equations with constant coefficient is considered with a
new approachment by using delta exponential function. Analysis of the method
is given in detailed and the advantage of this new approachment is to enable
us to investigate the solution of difference equations in the closed form
given in $\left( \ref{i2}\right) ,$ otherwise it has to be expanded by using
binomial expansion of the difference operator%
\begin{equation*}
\Delta ^{N}x\left( t\right) =\sum_{k=0}^{N}\left( -1\right) ^{k}\left( 
\begin{array}{c}
N \\ 
k%
\end{array}%
\right) x\left( t+N-k\right)
\end{equation*}%
and the method is supported with two difference eigenvalue problems; the
second-order Sturm-Liouville $\left( \ref{a2}\right) $ and the fourth-order
relaxation difference equations $\left( \ref{t1}\right) $. We find the sum
representation of solution of Sturm-Liouville difference problem and also,
we find the analytical solution of the fourth-order relaxation difference
problem.

Moreover, behaviors of eigenfunctions for the problems $\left( \ref{a2}%
\right) -\left( \ref{a3}\right) $ and $\left( \ref{a4}\right) -\left( \ref%
{a3}\right) $ are analyzed and illustrated by graphics and tables.. Firstly,
we show the behaviors of eigenfunctions while $q\left( n\right) =\frac{1}{%
\sqrt{n}}$ and eigenvalues are continuous in $Fig1$ and $Fig2,$ and we
observe that eigenfunctions are continous according to the eigenvalues. Also
we compare the behaviors of eigenfunctions in $Fig3.$ Then, we analyze the
behaviors of eigenfunctions for the specific eigenvalues $\lambda =1$ and $%
\lambda =2-\sqrt{2}$ in $Table1,$ $Table2$ and $Fig4$ and $Fig5$ and we
observe that eigenfunctions are discrete. We analyze the similar properties
for the problem $\left( \ref{t1}\right) -\left( \ref{t2}\right) $ in $Fig6$
while $q\left( n\right) =n,$ in $Fig7$ and $Table3$ while $\lambda =0.0625,$
in $Fig8$ and $Table4$ while $\lambda =0.1296$.

\end{document}